\documentclass[12pt]{amsart}
\usepackage{fullpage, amsmath, amscd, amssymb, amssymb}
\usepackage[all]{xy}

\newtheorem{The}{Theorem}
\newtheorem{Cor}[The]{Corollary}
\newtheorem{Def}{Definition}[section]
\newtheorem{Lem}[Def]{Lemma}
\newtheorem{Pro}[Def]{Proposition}

\newtheorem{Cla}[Def]{Claim}

\newtheorem{Rem}[Def]{Remark}
\newtheorem{Prob}{Problem}

\theoremstyle{remark}
\newtheorem{remark}{Remark}

\newcommand {\La}{\Lambda}
\newcommand {\la}{\lambda}
\newcommand {\ka}{\kappa}

\newcommand {\m}{\mathfrak{p}}

\newcommand {\len}{\ell}
\newcommand {\gln}{G_{\len^n}}
\newcommand {\kk}{k_F}

\newcommand {\Ok}{\mathfrak{o}_\len}

\newcommand {\Fb}{F}

\newcommand {\gb}{\mathbf{g}}
\newcommand {\cb}{\mathbf{c}}
\newcommand {\eb}{\mathbf{e}}

\newcommand {\ta}{\tau}
\newcommand {\f}{\phi}
\def \e{\epsilon}

\newcommand {\Db}{\mathbf{D}}
\newcommand {\Tb}{T}

\newcommand {\N}{\mathbb{N}}

\newcommand {\C}{\mathbb{C}}

\newcommand {\HH}{\mathcal{H}}

\newcommand {\LL}{\mathcal{L}}

\newcommand {\KK}{\mathcal{K}}
\newcommand {\Scal}{\mathcal{S}}

\newcommand {\FF}{\mathcal{F}}

\newcommand {\OO}{\mathfrak{o}}

\newcommand {\UU}{\mathcal{U}}

\newcommand {\dar}{\downarrow}

\newcommand {\rar}{\rightarrow}

\newcommand {\FFf}{\mathcal{F_{\phi}}}
\newcommand {\HHf}{\mathcal{H_{\phi}}}

\newcommand {\Span}{{\mathrm{Span}}}

\newcommand{\Gr}{\mathrm{Gr}}

\newcommand{\GL}{\mathrm{GL}}
\newcommand{\PGL}{\GL_n}

\newcommand{\End}{\mathrm{End}}
\newcommand{\Hom}{\mathrm{Hom}}
\newcommand {\coker}{\mathrm{coker}}

\newcommand{\ad}[2]{\Tb_{#2 \prec #1}}

\newcommand{\au}[2]{\Tb_{#2 \succ #1}}

\newcommand{\cont}[3]{[#1 \!\prec \!\dot{#2}]_{#3}}

\newcommand{\avoid}[3]{[#1 \!\prec \!\dot{#2} \!\pitchfork\! #3]}

\newcommand{\kova}[3]{\Tb_{#1 \hookrightarrow #2
\hookleftarrow #3}}

\newcommand {\rk}{\mathrm{rk}}

\newcommand{\justavoid}[2]{[\dot{#1} \!\pitchfork\! #2]}
\newcommand{\Ker}{\mathrm{Ker}}

\title{On some Geometric Representations of $\mathbf{GL_n(\OO)}$}
\author{Uri Bader and Uri Onn}
\thanks{Both authors were supported by the Israel Science Foundation, ISF grant no.
100146. The second author was partially supported by NWO grant no.
613.006.573} 


\begin{document}

\maketitle

\begin{abstract}
We study a family of complex representations of the group
$\GL_n(\OO)$, where $\OO$ is the ring of integers of a
non-archimedean local field $\Fb$. These representations occur in
the restriction of the Grassmann representation of $\GL_n(\Fb)$ to
its maximal compact subgroup $\GL_n(\OO)$. We compute explicitly
the transition matrix between a geometric basis of the Hecke algebra associated with the
representation and an algebraic basis which consists of its minimal idempotents. The transition matrix involves combinatorial invariants of
lattices of submodules of finite $\OO$-modules. The idempotents
are $p$-adic analogs of the multivariable Jacobi polynomials.
\end{abstract}



\section{Introduction}\label{introduction}

\subsection{Outline} Let $\Fb$ be a non-archimedean local field and $\OO$
its ring of integers. Let $\kk$ be the residue field and $\m$ the
maximal ideal of $\OO$. For $\len \in \N$, let $\Ok$ denote the
finite quotient $\OO/\m^\len$. The group $\PGL(\Fb)$ acts
transitively on $\Gr(m,n,\Fb)$, the Grassmannian of
$m$-dimensional subspaces of a fixed $n$-dimensional space, giving
rise to a complex representation of $\PGL(\Fb)$ on
$L^2(\Gr(m,n,\Fb))$. This work is about the restriction of this
representation to the maximal compact subgroup $\PGL(\OO)$. Since
the irreducible constituents of this representation are contained
in $\Scal(\Gr(m,n,\Fb)) \subset L^2(\Gr(m,n,\Fb))$, the dense
subspace of locally constant functions, we focus on the latter and
call it the \emph{Grassmann representation} of $\GL_n(\OO)$. The
Grassmann representation has a multiplicity free decomposition to
irreducible representations
\begin{equation}\label{decomposition}
\Scal(\Gr(m,n,\Fb)) = \bigoplus_{\lambda \in \Lambda_m}
\UU_\lambda^F,
\end{equation}
where $\Lambda_m$ stands for partitions of at most $m$ parts, see
\cite{p1}. Each irreducible constituent $\UU_\lambda^F$ contains a
unique (normalized) $P_x$-spherical vector $e_\lambda^F$, where $P_x$ is a stabilizer of a point $x \in \Gr(m,n,\Fb)$.
These functions are the non-archimedean analogs of
the multivariable Jacobi polynomials which arise in the analogous
setup when $F$ is either $\mathbb{R}$ or $\C$, and will therefore
be called the $p$-adic multivariable Jacobi functions.
Algebraically, after the appropriate normalization, they form a
basis  of the Hecke algebra
$\HH_m=\End_{\GL_n(\OO)}\bigl(\Scal(\Gr(m,n,F))\bigr)$ which consists of minimal idempotents. The algebra
$\HH_m$ can be identified with the convolution algebra $\Scal(P_x
\backslash \GL_n(\OO)/P_x)$ of bi-$P_x$-invariant locally constant
functions on $\GL_n(\OO)$. The latter description comes with a natural
geometric basis: characteristic functions of the double cosets.
The main result in this paper is an explicit computation of the
transition matrix between these geometric and algebraic bases
(Theorem \ref{main.theorem}).

The Grassmann representation is filtered by the finite dimensional
$I_\len$-invariant subspaces
\begin{equation}\label{filtration.grass}
(0) \subset \Scal(\Gr(m,n,\Fb))^{I_1} \subset \cdots \subset
\Scal(\Gr(m,n,\Fb))^{I_\len} \subset \cdots \subset \Scal(\Gr(m,n,\Fb)),
\end{equation}
where $I_\len=\Ker\{\GL_n(\OO) \to \GL_n(\OO_\len)\}$. In fact,
$\Scal(\Gr(m,n,\Fb))=\varinjlim \Scal(\Gr(m,n,\Fb))^{I_\len}$,
hence the problem can be translated into a finite problem:
analysis of the representation of the finite group
$\GL_n(\OO_\len)$ in $\Scal(\Gr(m,n,\Fb))^{I_\len}$. The latter
can in turn be identified with $\C(\Gr(m,n,\OO_\len))$, the space
of complex valued functions on the (finite) Grassmannian of free
$\OO_\len$-submodules of rank $m$ in $\OO_\len^n$, with its
natural $\GL_n(\OO_\len)$-action.

\subsection{Context of the problem} To put things into perspective,
we briefly describe the archimedean \cite{james-consta} and quantum \cite{Stokman2} counterparts of
the Grassmann representation, see also \cite{onn1,OS} for more details. These are representations of the
orthogonal group $\mathrm{O}_n$ arising from its action on
$\Gr(m,n,\mathbb{R})$, of the unitary group $\mathrm{U}_n$
arising from its action on $\Gr(m,n,\mathbb{C})$, and of the
quantized enveloping algebra $\mathrm{U}_q(\mathfrak{gl}_n)$ on
the quantum Grassmanian \cite{Stokman2}. They all have similar decompositions to \eqref{decomposition}, indexed
by the same set $\Lambda_m$, giving rise to zonal
spherical functions $e_\lambda^\mathbb{R}$ and $e_\lambda^\C$, which are special cases of the multivariable Jacobi polynomials, and
$e_\lambda^q$ which are known as the multivarable little $q$-Jacobi polynomials \cite{Stokman1}.

In \cite{Stokman1,KS} it is shown that by taking appropriate
limits, the multivariable $q$-Jacobi polynomials $e_\lambda^q$
degenerate to the multivariable Jacobi polynomials which
specialize to $e_\lambda^\mathbb{R}$ and $e_\lambda^\C$, and in
\cite{onn1} it is further shown that they degenerate to the
$p$-adic multivariable Jacobi functions $e_\la^F$ which are
studied in the present paper (see also \cite{Haran} and \cite{KO}
for the projective space case). Besides the aesthetic nature of
these limits, in which the quantum zonal spherical functions
degenerate to the zonal spherical functions over all local fields,
archemedeans and non-archimedeans, they have been used in
\cite{OS} to compute the dimensions of the irreducible
constituents in \eqref{decomposition}, and no other direct
computation is known at present.

In a different direction, the multivariable $p$-adic Jacobi
functions generalize the $q$-Hahn polynomials for
$q=|\OO/\m|=|\kk|$ (see e.g. \cite{GR} for the precise
definition), which under the appropriate interpretation coincide
with $e^F_{\la}$, for $\la=(1^j)$ with $0 \le j \le m$. These
functions form a basis of $\End_{\GL_n(\kk)}\C(\Gr(m,n,\kk))$,
which captures the first term in the filtration
\eqref{filtration.grass}. See \cite{delsarte-H,dunkl} for more
details on this special case.

The irreducible representations $\UU_\la^F$ are studied in \cite{p1} in great detail. Their precise identification involves the study a wider
family of geometric representations of the group $\GL_n(\OO_\len)$ which arise from its action on $\Gr(\la,\OO_\len^n)$, the Grassmannian of submodules
of $\OO_\len^n$ of type $\la$. Isomorphism types of submodules of $\OO_\len^n$ are classified by partitions, which is the underlying reason behind the appearance of $\Lambda_m$ in
\eqref{decomposition}.

\subsection{Content of the paper} Section \ref{hecke-alg} is devoted to representations and Hecke algebra which arise from the action
of the finite quotients $\GL_n(\OO_\len)$ on Grassmannians of submodules of $\OO_\len^n$. The main tools, i.e. geometrically defined intertwining operators, are
 described and developed. Most of this section is an adaptation of relevant results and ideas from \cite{p1} which is the foundational background for this work.

  Sections \ref{sec:c-g} and \ref{sec:c-e} are devoted to transition matrices between various bases of the Hecke algebras. The main tools which are used are
   combinatorial invariants of the lattice of submodules in free $\OO_\len$-modules, the Euler characteristic of the simplicial complex associated with flags of
   submodules, and explicit computations with $q$-binomial coefficients.

 In section \ref{section:infinite} the finite results are transferred to $\Scal(\Gr(m,n,\Fb))$ and its Hecke algebra $\HH_m$. Special attention is given to
 the limiting process from the algebraic and topological aspects.

   Section \ref{problems} is devoted to related topics and open problems. In the appendix several claims on modules over discrete valuation rings are proved, which
   we suspect to be known, but could not find an adequate reference.

\subsection{Acknowledgements}
We are most grateful to Amos Nevo for hosting and encouraging
this research. We warmly thank Shai Haran for stimulating
discussions, his interest and encouragement.


\section{Hecke algebras associated to finite Grassmannians}\label{hecke-alg}

Let $\OO_\len=\OO/\m^\len$ and let $G_{\len^n}$ denote the automorphism group of a free $\OO_\len$-module of rank $n$. Upon a choice of a basis
$G_{\len^n}$ can be identified with $\GL_n(\OO_\len)$. Recall that any finitely generated $\OO_\len$-module is isomorphic to $\OO_\la=\oplus_{i=1}^{r}
\OO/\m^{\la_i}$ for some $\la=(\la_i)$ where $\len \ge \la_1\ge\cdots\ge\la_r \ge 0$ since $\OO$ is a discrete valuation ring. We call $\la$ the (isomorphism) {\em type} of that module. The {\em length} of the partition $\la$, i.e. the number of its nonzero parts, is the rank of the module $\OO_\la$, and the {\em height} of a partition is its largest part. The set of partitions of length at most $m$ is denoted $\La_m$.  The set of types is equipped with a
natural partial order: $\mu \le \nu$ whenever a module of type
$\mu$ can be embedded in a module of type $\nu$. In terms of
the corresponding Young diagrams it amounts to inclusion of the corresponding
diagrams. Here and in the sequel we follow the notations of \cite[II.1]{macdonald}.

In this section we analyze representations and Hecke algebras which arise from Grassmannians of submodules of the free $\OO_\len$-module of rank $n$.
On one hand they generalize the case $\len=1$ which is studied in \cite{dunkl}, and on the other hand they form the
crucial ingredient for understanding the Grassmann representation $\Scal\left(\Gr(m,n,\Fb)\right)$ of $\GL_n(\OO)$.

\subsection{Definition and realization of the Hecke algebra}

Let $\LL(\OO_\len^n)$ denote the lattice of submodules of $\Ok^n$. The
group $G_{\len^n}$ acts on the lattice $\LL(\OO_\len^n)$ which is a disjoint union of the $G_{\len^n}$-invariant subsets
\[
X_{\la}=\Gr(\la,\OO_\len^n)=\{x \in \LL(\OO_\len^n) ~|~ x \simeq \OO_\la\} \qquad (\la \in \La^\len_n),
\]
where
$\La^\len_n=\{\la \in \La_n~|~\text{height}(\la) \le \len\}$
stands for isomorphism types of elements in $\LL(\OO_\len^n)$. Let $\tau:\LL(\OO_\len^n) \rightarrow \La^\len_n$ be the
{\em type map} which assigns to each module its isomorphism type.

For each $\la \in \La_n^\len$, let $\FF_{\la}=\C(X_{\la})$ be the
complex permutation representation of $G_{\len^n}$ arising from
its action on $X_{\la}$. We equip $\FF_{\la}$ with the inner
product associated to the counting measure on $X_{\la}$. The
following claim, which will be proved in \S\ref{section:infinite},
highlights the relevance of the representations $\FF_\la$ with
$\la=(\len^m)=(\len,\ldots,\len)$.
\begin{Cla}\label{isom.grass.and.limit}  There exist an isomorphism of $\GL_n(\OO)$-representations
\[
\Scal\left(\Gr(m,n,\Fb)\right) \simeq \varinjlim \FF_{\len^m}.
\]
\end{Cla}
Apart from their role in Claim \ref{isom.grass.and.limit}, the representations $\FF_{\len^m}$ are distinguished among all $\FF_{\la}$ since
\begin{enumerate}
\item [(a)] They are multiplicity free, and
\item [(b)] The number of their irreducible constituents is $|\La_m^\len|=\big({\len+m \atop m}\big)$, in particular, it is independent of $\OO$.
\end{enumerate}
This is proved in \cite{p1} in greater generality, but in order to be as self contained as possible we shall explain it in detail. Let
\[
\HH_{\len^m}=\End_{\gln}(\FF_{\len^m})
\]
stand for the (Hecke)
algebra of $\gln$-invariant endomorphisms of $\FF_{\len^m}$. Assertions (a) and (b) would follow once we show that the algebra $\HH_{\len^m}$ is isomorphic
to the algebra $\C^{\La_m^\len}$ of complex valued functions on $\La_m^\len$ with pointwise multiplication. To prove that, we look at an alternative
description of $\HH_{\len^m}$ which is of geometric flavor. The algebra $\End_\C(\FF_{\len^m})$ can be identified with $\C(X_{\len^m} \times X_{\len^m})$ by interpreting
any function $f:X_{\len^m} \times X_{\len^m} \to \C$ as a summation kernel $T_f:\FF_{\len^m} \to \FF_{\len^m}$, that is, $T_f(h)(x)=\sum_{y \in X_{\len^m}}f(x,y)h(y)$. The map $f \mapsto T_f$ is just the identification of matrices with endomorphisms, which is also equivariant with respect to the natural $\gln$-action on both algebras:
$[g\cdot T](h)=g(T(g^{-1}h))$ and $[g \cdot f](x,y)=f(gx,gy)$ for $g \in \gln$, $T \in \FF_{\len^m}$, and $f \in \C(X_{\len^m} \times X_{\len^m})$. Taking $\gln$-invariants gives
\begin{equation}\label{identification.hecke}
\HH_{\len^m} \simeq \C(X_{\len^m} \times_{\gln} X_{\len^m}).
\end{equation}

\begin{Cla}\label{ident.cosets.types} For $m \le n/2$ there exist a bijection
\begin{equation*}
 \begin{split}
 X_{\len^m} \times_{\gln} X_{\len^m}  &\overset{\sim}{\longrightarrow} \La_m^\len \\
   \gln(x,y) &\longmapsto   \tau(x \cap y),
\end{split}
\end{equation*}
\end{Cla}

\begin{proof} The fact that $\gln$ preserves the module
structure implies that this map is well defined. It is onto due
to the assumption $m \le n/2$ which gives enough room to realize
any type $\la$ as intersection of two $\OO_\len$-modules of type
$\len^m$. It is one-to-one because any abstract isomorphism
between $x \cap y$ and $x' \cap y'$ can be lifted to an element
$g \in \gln$ such that $(x',y')=(gx,gy)$, using
%
\begin{Lem}\label{symmetricity} Let $z \subset E$ and $z' \subset E'$ be modules such that $z \simeq z'$ and $E \simeq E' \simeq \OO_\len^j$. Then any isomorphism of $\OO_\len$-modules $z \to z'$ can be extended to an isomorphism $E \to E'$.
\end{Lem}
\begin{proof} The ring $\OO_\len$ is self injective\footnote{This is well known, but can also be easily verified using Baer's criterion which reduces the injectivity verification to a trivial calculation.}, therefore, the module $\OO_\len^j$ is injective, and hence the embedding $z \overset{\sim}{\to} z' \hookrightarrow E'$ can be extended to a map $E \to E'$. It is easy to see that among such extensions exist one-to-one extensions
which must be surjective as well due to the finiteness of $\OO_\len$.
\end{proof}
Going back to the argument above, the isomorphism $x \cap y \simeq x' \cap y'$ can be (simultaneously) extended to isomorphisms $x \simeq x'$ and $y \simeq y'$ by using the lemma for $j=m$. These two isomorphisms glue to an isomorphism $x+y \simeq x'+y'$, which by using the lemma once more for $j=n$, proves the claim.
\end{proof}
\begin{Cor}\label{km} The algebra $\HH_{\len^m}$ is semisimple, commutative and of dimension $|\La_m^\len|$, hence
\begin{equation}\label{identification.algebraic}
\HH_{\len^m} \simeq \C^{\La_m^\len}.
\end{equation}
\end{Cor}
\begin{proof}  Given the bijection of Claim \ref{ident.cosets.types} we get that both $\HH_{\len^m}$ and $\C^{\La_m^\len}$ have the same dimension:
$|\La_m^\len|=|\{ \text{isomorphism types of submodules of $\OO_\len^m$}\}|$. Both are semisimple, thus the only non obvious issue is the commutativity of the algebra $\HH_{\len^m}$, which
follows from Gelfand's trick. That is, we have that $(x,y)$ and $(y,x)$ are in the same $\gln$-orbit since $\tau(x \cap y)=\tau(y \cap x)$, hence the identity is an anti-isomorphism of the algebra.
\end{proof}

It will be convenient to denote from now on
$\phi=\len^m$ and $\Phi=\len^n$. To make the link with the terminology of \cite{p1}, note that Lemma \ref{symmetricity} shows that rectangular types such as $\phi$ and $\Phi$
are symmetric (Definition 2.1 in {\em loc. cit.}), and Claim \ref{ident.cosets.types} shows that $(\phi,\Phi)$ form a symmetric couple (Definition 2.2 in {\em loc. cit.}), provided that $m \le n/2$, which is our assumption throughout.

Using the isomorphisms \eqref{identification.hecke} and \eqref{identification.algebraic}, we obtain two natural bases for $\HH_{\phi}$.
The first, which comes from the r.h.s of \eqref{identification.hecke}, consists of the operators corresponding to
 characteristic functions of the orbits $\{\gb_\la~|~ \la \le \phi \}$, and will be call the {\em geometric basis}, and the second, which comes from the r.h.s of \eqref{identification.algebraic}, consists of idempotents $\{\eb_\la~|~ \la \le \phi \}$, and will be called the {\em algebraic basis}. The precise meaning of the indexing of the algebraic basis will follow from Theorem \ref{ideals-reps} below. In order to connect these bases we need a fine analysis of intertwining operators which is the theme of the next subsection.

\subsection{Geometric intertwiners and bases for the Hecke algebra}
Define the following operators
\begin{enumerate}
\item [(a)] For each
pair of types $\la \le \mu$ let
\begin{equation*} \au{\la}{\mu}: \FF_{\la} \rightarrow \FF_{\mu}, \qquad \au{\la}{\mu}h(y)=\sum_{ x \subset
 y } h(x)  \qquad (y \in X_{\mu})
\end{equation*}
\begin{equation*}
 \ad{\mu}{\la}: \FF_{\mu} \rightarrow \FF_{\la}, \qquad \ad{\mu}{\la}h (x)=\sum_{y \supset x } h(y)  \qquad (x \in X_{\la}).
\end{equation*}

\medskip

\item [(b)] For types $\la,\mu \le
\nu$, let $\kova{\la}{\nu}{\mu}:\FF_{\mu} \to \FF_{\la}$ be the operator
\begin{equation*}
\kova{\la}{\nu}{\mu}h(x)=\sum_{\{y \in X_\mu ~|~ \tau(y+x)=\nu \}}h(y) \qquad (x \in X_{\la}).
\end{equation*}

\medskip

\item [(c)] The aforementioned operators $\gb_{\la} \in \HH_\f$ ($\la \le \f$) are explicitly defined by
\[
\gb_{\la}h(x)=\sum_{\{y ~|~ \tau(y \cap x)=\la\}}h(y)\qquad (x \in X_{\f}).
\]

\end{enumerate}

Note that all these operators commute with the $G_\Phi$-action, and that
$\au{\la}{\mu}$ and $\ad{\mu}{\la}$ form an adjoint pair. In order to minimize confusion, we follow the rule that whenever
an operator is labeled with a diagram (e.g. $\kova{\la}{\nu}{\mu}$), it acts from the space
indexed by the \emph{right} type of the diagram ($\FF_\mu$) to the space
indexed by the \emph{left} type ($\FF_\la$).

\medskip

In order to find the transition matrix between the geometric basis $\{\gb_\la\}$ and the algebraic basis $\{\eb_\la \}$ of $\HH_\f$, we introduce
a third basis which is defined by $\cb_{\la}=\au{\la}{\f}\ad{\f}{\la}$  ($\la \le \f$), and will be called the {\em cellular basis}.
The geometric basis can be viewed as averaging operators along \lq spheres\rq, whereas the cellular basis can be viewed as weighted averaging operators on \lq balls\rq. More specifically, the following upper triangular relation holds \cite[\S3.4.2]{p1}
\begin{equation}\tag{\textbf{c-g}}
\cb_{\la}=\sum_{\f \ge \nu \ge \la} \Bigl({\nu \atop
\la}\Bigr)\gb_{\nu},
\end{equation}
which in turn proves that $\{\cb_{\la}\}$ is indeed a basis. Here $\bigl({\nu
\atop \la}\bigr)$ is the number of submodules of type $\la$
contained in a module of type $\nu$. For each $\la \le \f$ set
\begin{equation*}
\HH_{\f}^{\la}=\Span_{\C}\{\cb_{\la'} ~|~ \la' \le \la \}, \quad \text{and}
\quad \HH_{\f}^{\la^{-}}=\Span_{\C}\{\cb_{\la'} ~|~ \la' < \la \}.
\end{equation*}
The following is proved in \cite{hill}, and in a greater
generality in \cite{p1}.
\begin{The} \label{ideals-reps} $\HH_{\f}^{\la}$ and $\HH_{\f}^{\la^-}$ are ideals $\forall
 \la \le \f$, hence, $\bigl\{\KK_{\la}=\HH_{\f}^{\la}/\HH_{\f}^{\la^-}\bigr\}_{\la \le
\f}$ is a complete set of irreducible representations of $\HHf$.
\end{The}
\begin{Cor}\label{triangularity.c-e} If $\eb_{\la}$ is the idempotent in $\HHf$
corresponding to $\KK_{\la}$ for all $\la \le \f$, then there
exist a lower triangular matrix $(A_{\la\ka})$ such that
\begin{equation}\tag{\textbf{c-e}}
\cb_{\la}=\sum_{\ka \le \la} A_{\la \ka}\eb_{\ka}.
\end{equation}
\end{Cor}

The cellular basis appears as a bridge between the geometric and
algebraic bases. It is upper triangular with respect to the former
and lower triangular with respect to the latter. In the next
subsection we shall use it to compute the idempotents explicitly.

\begin{remark} Theorem \ref{ideals-reps} can be used to label the
irreducible representations of $\HHf$. The one-dimensional
$\HHf$-module $\KK_{\la}$ is the unique $\HHf$-module which is
annihilated by all $\{\HH_{\f}^{\mu}~|~\mu < \la\}$ and not
annihilated by $\HH_{\f}^{\la}$. In view of the well known
dictionary between representations of the group which occur in
$\FFf$ and modules of the Hecke algebra $\HHf$, we can now label
the irreducibles in $\FFf$ by: $\UU_{\la} \leftrightarrow
\KK_{\la}$. Moreover, by the definition of $\cb_{\la}$ as the
composition $\au{\la}{\f} \ad{\f}{\la}$, the annihilation
criterion above translates to the fact that $\UU_{\la}$ occurs in
$\FF_{\la}$ and does not occur in $\FF_{\mu}$ for $\mu < \la$.
\end{remark}

\section{Transition matrix: cellular to geometric}\label{sec:c-g}

In this section we invert the relation \textbf{(c-g)} and compute it explicitly. It consists of two parts, an abstract inversion using properties
of the lattice of submodules and an explicit calculation in terms of $q$-binomial coefficients, where $q$ is the cardinality of the residue field $\kk$.

\subsection{An abstract inversion} Two $\OO$-module monomorphisms $i:x \hookrightarrow y$ and $i':x' \hookrightarrow y'$ are said to be equivalent if there are
isomorphisms $x \simeq x'$ and $y \simeq y'$ such that the
following diagram is commutative
\begin{equation}
\begin{matrix}
x & \underset{i}{\hookrightarrow} & y \\
\downarrow & & \downarrow  \\
x' & \underset{i'}{\hookrightarrow} & y'
\end{matrix}
\end{equation}
Assuming the isomorphism types of $x$ and $y$ are $\lambda$ and
$\nu$ correspondingly, we denote the equivalence class of $i:x
\hookrightarrow y$ by $i: \lambda \hookrightarrow \nu$, and let $\bigl({\nu=\nu
     \atop i:\la \hookrightarrow \nu}
     \bigr)$ stand for the number of submodules of type $\la$ in a module of type $\nu$ with embedding type $i$.

Let $y$ be a finite module over $\OO$. Denote by $\LL(y)$ the
lattice of submodules of $y$. One naturally associates a
simplicial complex to $y$, denoted $C(y)$, with vertices
consisting of the non-trivial submodules of $y$ (all but $0$ and
$y$). The simplices of $C(y)$ are given by flags
\[ \{(y_1,y_2,\ldots,y_m)~|~0\subset y_1\subset y_2\subset\cdots
   \subset y_m\subset y \}
\]
We denote the Euler characteristic of $C(y)$ by $\chi(y)$.

For inverting the relation (\textbf{c-g}) we note that coefficients $\left({\nu \atop
   \la}\right)$ coincide with $\hat{\zeta}(\la,\nu)$ in the notation of \cite[\S 2.2]{p1} (Proposition 2.5 and the discussion afterwards). Its inverse (denoted $\hat{\mu}(\la,\nu)$) is given by
  \[ \sum_{\la\underset{i}{\hookrightarrow}\nu~}
     \Bigl({\nu=\nu
     \atop \la\underset{i}{\hookrightarrow}\nu}
     \Bigr) \chi(\coker(i)), \qquad \text{\cite[Claim 2.6]{p1}}.
  \]
We get
\begin{equation}\label{inverted.g-c}
 \gb_{\nu}=\sum_{\la\le\nu}\sum_{i}
     \Bigl({\nu=\nu
     \atop \la\underset{i}{\hookrightarrow}\nu}
     \Bigr) \chi(\coker(i))\cb_{\la}.
 \end{equation}
The following Lemma shows that many of the terms in \eqref{inverted.g-c} vanish.

\begin{Lem} \label{euler} If $\m y\neq (0)$ then $C(y)$ is contractible, in
particular, $\chi(y)=0$.
\end{Lem}

\begin{proof}
Denote $C=C(y)$, and $D\subseteq C$ the subcomplex spanned by the
subset of vertices $\{x~|~\m y \subseteq x\subset y\} \subseteq
C_0$. $D$ is a cone over the vertex $\m y$, hence contractible. The
function
\[ \varphi:\LL(y)\rightarrow \LL(y), \quad \varphi(x)=x+\m y \]
extends to a retraction $\varphi_*:C\rightarrow D$. The function
\[ \Psi : C\times [0,1] \rightarrow C \qquad
   \Psi(c,t)=(1-t)c+t\varphi_*(c)
\]
establishes a deformation retract from $C$ to $D$. Therefore $C$
is contractible.
\end{proof}

The value of $\chi(y)$ depends only on the isomorphism type of $y$ hence we shall use the notation $\chi(\la)$ for $\la \in \Lambda_n$. The vanishing
of $\chi(\la)$ when $\OO_\la$ is not a vector space will be written in short as $\chi(\la)=0$ if $\m\la=0$.

\subsection{Explicit calculation}\label{explicit.1}

It will be useful to use another set of coordinates for
elements in $\La^\len_n$, obtained by transposed diagrams
$\la'=(\la'_j)$, defined by $\la'_j=|\{i:\la_i \ge j\}|$. The module-theoretic interpretation of the $\la_j'$'s is given in
\cite[II.1(1.4)]{macdonald}. For every partition $\xi$, let $n(\xi)=\sum (i-1)\xi_i$ and $|\xi|=\sum \xi_i$. Let $q=\big|\OO/\m\big|=|\kk|$, and for $m,n \in \N$ set
\[
[n]_q=1-q^{-n} , \qquad \qquad [0]_q=1, \qquad \qquad \quad
\]
\[
\qquad \qquad  [n]_q!=[n]_q[n-1]_q \cdots [1]_q , \qquad \Bigl[{m
\atop n}\Bigr]_q=\frac{[n]_q!}{[m]_q![n-m]_q!}.
\]
The subscript $q$ will be occasionally omitted from the notation.
\begin{Pro}
\[ \gb_{\la}=\sum_{\{\nu|\f \ge \nu \ge \la \ge \m\nu\}}
   (-1)^{|\nu|-|\la|}q^{n(\nu)-n(\la)} \prod_{i \ge 1}
   \biggl[{\nu_i'-\nu_{i+1}' \atop \nu_i'-\la_i'}\biggr] \cdot
   \cb_{\nu},
\]
\end{Pro}

\begin{proof}
Lemma~\ref{euler} combined with a
well known formula \cite{foundation-i,chen-rota} giving the Euler characteristic of the Tits
building associated to $\GL_n(\kk)$, gives for every module type
$\la$,
\begin{equation} \label{chi}
 \chi(\la)=\left\{
    \begin{array}{lr}
    0 & \quad \m\la\neq 0 \\
    (-1)^{\dim_{\kk}(\la)} q^{\bigl({\dim_{\kk}(\la) \atop 2}\bigr)} & \quad \m\la= 0
    \end{array}\right.
\end{equation}

whereas the number
\[ \sum_{\{\la\underset{i}{\hookrightarrow}\nu~|~\m\coker(i)=0\}}
     \Bigl({\nu=\nu
     \atop \la\underset{i}{\hookrightarrow}\nu}
     \Bigr)
\]
is exactly the Hall coefficient $G^{\nu}_{\la,1^{(|\nu|-|\la|)}}$.
The latter is explicitly computed in \cite[II.4]{macdonald}.
\end{proof}


\section{Transition matrix: cellular to idempotents}\label{sec:c-e}

\subsection{Abstract description} Let $E$ be a fixed $\OO$-module of type $\Phi=\len^n$. We say
that two submodules are \emph{transversal} if their intersection
is zero. Let $\ka$, $\la$ and $\mu$ be types of modules and let
$x_{\mu}$ and $x_{\ka}$ be two transversal submodules of $E$ of
types $\mu$ and $\ka$ respectively. Let

\begin{itemize}
\item $[\mu \prec \dot{\la} \pitchfork \ka]_{\Phi}$ ~~be the number
of submodules of type $\la$ in $F$ which contain $x_{\mu}$ and are
transversal to $x_{\ka}$.

\item $[\mu \prec \dot{\la}]_{\Phi}$ ~~be the number of submodules of
$F$ of type $\la$ which contain a given submodule of type $\mu$
(in the above notation this is $[\mu \prec \dot{\la} \pitchfork 0
]_{\Phi}$).

\end{itemize}
Note that both $[\mu \prec \dot{\la} \pitchfork \ka]_{\Phi}$ and $[\mu \prec \dot{\la}]_{\Phi}$ are well defined by Lemma \ref{symmetricity}.

\medskip

Since the cellular structure agrees with the idempotent
decomposition, we already know by Corollary
\ref{triangularity.c-e} that there exist a lower triangular matrix
$A_{\la\ka}$ such that the relation (\textbf{c-e}) above holds. We
have already seen that the transition matrix (\textbf{c-g}) from
the geometric basis to the cellular basis depends only on
geometric invariants of the lattice of submodules in a very simple
way. This is also the case for the cellular-idempotents transition
matrix.

\begin{The}\label{cell-idem} $A_{\la\ka} = \cont{\ka}{\la}{\f} \avoid{\la}{\f}{\ka}_\Phi$.
\end{The}
Note that Theorem 4 is a generalization of \cite[Theorem 6]{p1} in
which the same statement is proved for $\lambda$'s with
rectangular shape. Our strategy is to analyze the multiplication
in the algebra with respect to the cellular basis. Let
$B_{\la\mu}^{\nu}$ be multiplication coefficients with respect to
the cellular basis
\begin{equation}\label{mtable}
 \cb_{\la} \cdot \cb_{\ka} = \sum_{\nu \le \la \wedge \ka}
 B_{\la\ka}^{\nu} \cb_{\nu}.
 \end{equation}
Observe that $B^{\ka}_{\la\ka}=A_{\la\ka}$ for $\ka \le \la$. The
following lemma follows immediately from \cite[Lemma
3.6]{p1}.
\begin{Lem}\label{composition} \qquad
\begin{enumerate}
    \item ~$\au{\ka}{\f}\au{\nu}{\ka}= \cont{\nu}{\ka}{\f}
    \au{\nu}{\f}$.
    \item ~$\ad{\ka}{\nu}\ad{\f}{\ka}=\cont{\nu}{\ka}{\f}\ad{\f}{\nu}$.
\end{enumerate}
\end{Lem}

\noindent Substituting $\cb_{\eta}= \au{\eta}{\f}\ad{\f}{\eta}$ in
equation \eqref{mtable} and using Lemma \ref{composition} gives
(assume $\ka \le \la$):
\begin{equation}\label{ddd}\begin{split}
 \au{\la}{\f}\Bigl(&\ad{\f}{\la}\au{\ka}{\f}\Bigr)\ad{\f}{\ka} \\
 &=\au{\la}{\f} \Bigl(\sum_{\nu \le \ka}
\frac{B_{\la\ka}^{\nu}}{\cont{\nu}{\la}{\f} \cont{\nu}{\ka}{\f}}
\au{\nu}{\la}\ad{\ka}{\nu}\Bigr)\ad{\f}{\ka}
\end{split}\end{equation}

Let $\FF_{\la}^{\bullet} = \text{Im}(\ad{\f}{\la})$. Observe that the set of maps
$\Delta=\Delta_{\la\ka}= \{\au{\nu}{\la}\ad{\ka}{\nu}\}_{\nu \le
\la \wedge \ka}$ when restricted to $\FF_{\ka}^{\bullet}$, forms
a basis for $G_\len$-maps $\FF_{\ka}^{\bullet} \rightarrow
\FF_{\la}^{\bullet}$. Indeed, these maps are independent after
being composed with $\au{\la}{\f}$ on the left and $\ad{\f}{\ka}$
on the right and using lemma \ref{composition}. We get the
following identity on $\FF_{\la}^{\bullet}$:
\begin{equation}\label{eee}\begin{split}
 \ad{\f}{\la}\au{\ka}{\f}&=\sum_{\nu \le \ka}
\frac{B_{\la\ka}^{\nu}}{\cont{\nu}{\la}{\f} \cont{\nu}{\ka}{\f}} \au{\nu}{\la}\ad{\ka}{\nu}\\
&= \frac{A_{\la\ka}}{\cont{\ka}{\la}{\f}} \au{\ka}{\la} + \{
\text{terms with} ~\nu < \ka \}
\end{split}\end{equation}
Denote the coefficient of an operator $S$ w.r.t to a basis element
$D \in \Delta$ by $\langle S, D\rangle_{\Delta}$.
\begin{Def}\label{good}
 A triple $\ka \le \la \le \eta$ is called \emph{good} if
 \[\langle \ad{\eta}{\la} \au{\ka}{\eta},
\au{\ka}{\la}\rangle_{\Delta} = \avoid{\la}{\eta}{\ka}_{\Phi}\]
\end{Def}
\noindent Combining Definition \ref{good} with \eqref{eee} we see
that Theorem \ref{cell-idem} is equivalent to:
\begin{The}\label{heart} $\ka \le \la \le \f$ is a good triple.
\end{The}
We would like to take a small pause and
explain the strategy which we undertake. The idea is to show that
it is enough to find a path connecting $\la$ and $\f$ in the
segment $[\la,\f]$ which can be paved with good triples, and then
exhibit such path. More precisely, we follow three steps:
\begin{enumerate}
    \item Given $\la_0 \le \la_1 \le \cdots \le \la_r \le \f$ such that
    $\ka \le \la_{i-1} \le \la_i$ is good for all $1 \le i \le r$ and also $\ka \le \la_r \le
    \f$ is good, implies that so is $\ka \le \la_0 \le \f$.

    \item $\ka \le \la_1 \le \la_2$ is good
    whenever $\la_2$ covers $\la_1$ (i.e. $[\la_1,\la_2]=\{\la_1,\la_2\}$) and has the same rank.

    \item $\ka \le \e \le \f$ is good whenever $\e$ and $\f$ are symmetric.
 \end{enumerate}
Note that a Jordan-H\"{o}lder sequence of types from $\la$ to a
symmetric type $\e$ of the same rank gives an appropriate path:
$\la = \la_0 \le \cdots \le \la_r=\e \le \f$.

\bigskip

\noindent \textbf{Step 1}

\begin{Lem} If $\ka \le \la \le \theta$ and $\ka \le \theta \le \f$ are good so is $\ka \le \la \le \f$.
\end{Lem}

\begin{proof}
\begin{equation*}\begin{split}
\langle \ad{\f}{\la}\au{\ka}{\f}, \au{\ka}{\la}\rangle_{\Delta} &=
\frac{1}{\cont{\la}{\theta}{\f}} \langle
(\ad{\theta}{\la}\ad{\f}{\theta})\au{\ka}{\f},
\au{\ka}{\la}\rangle_{\Delta} \\
&= \frac{\avoid{\theta}{\f}{\ka}}{\cont{\la}{\theta}{\f}} \langle
\ad{\theta}{\la}(\ad{\ka}{\theta}+{\text{lower terms})},
\au{\ka}{\la}\rangle_{\Delta} \\ &=
\frac{\avoid{\theta}{\f}{\ka}}{\cont{\la}{\theta}{\f}}\langle
\ad{\theta}{\la}\ad{\ka}{\theta}, \au{\ka}{\la}\rangle_{\Delta}
\\ &= \frac{\avoid{\theta}{\f}{\ka}
\cont{\la}{\theta}{\ka}}{\cont{\la}{\theta}{\f}} = \frac{[\la
\!\prec\! \overset{\centerdot}{\theta} \!\prec\!
\overset{\centerdot}{\f} \!\pitchfork\!
\mu]}{\cont{\la}{\theta}{\f}} = \avoid{\la}{\f}{\ka}_{\Phi}
\end{split}
\end{equation*}
\end{proof}
Using this lemma $r$ times completes the proof of step 1. Note
that the lemma can be used successively only from the 'top' to the
'bottom'.
\bigskip

\noindent \textbf{Step 2}

\noindent Let $\la_1$ and $\la_2$ be types of same rank and assume
that $\la_2$ covers $\la_1$. Let $\ka \le \la_1$. We want to show
that $\langle \ad{\la_2}{\la_1} \au{\ka}{\la_2},
\au{\ka}{\la_1}\rangle_{\Theta} = \avoid{\la_1}{\la_2}{\ka}$.
However, the assumption that $\la_1$ and $\la_2$ have the same
rank guarantees that any module of type $\la_2$ containing a
module of type $\la_1$ which is transversal w.r.t. a module of
type $\ka$ inherits this transversality. Hence, the requirement to
avoid $\ka$ is redundant and
$\avoid{\la_1}{\la_2}{\ka}=\cont{\la_1}{\la_2}{}$. We start by
expanding the product $\ad{\la_2}{\la_1} \au{\ka}{\la_2}$:
\begin{equation}\label{H}
\ad{\la_2}{\la_1} \au{\ka}{\la_2}=
a_{\la_1}\kova{\la_1}{\la_1}{\ka}+a_{\la_2}\kova{\la_1}{\la_2}{\ka}
\end{equation}
The assumptions on $\la_1$ and $\la_2$ assures that no other terms
appear in \eqref{H}. Evidently
$\au{\ka}{\la_1}=\kova{\la_1}{\la_1}{\ka}$ and
$a_{\la_1}=\cont{\la_1}{\la_2}{}$. We are therefore reduced to
showing that the $\Delta$ expansion of the second term in
\eqref{H} does not contain a multiple of $\au{\ka}{\la_1}$. This
is accomplished by claim \ref{WWW}. Let $\ka_1$ be the unique type
which can (possibly) complete a cartesian diagram (see the first
part of claim \ref{equivalence} in appendix A):
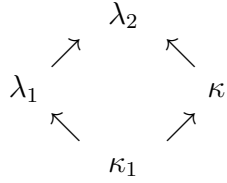
\begin{equation*}
\begin{matrix}   & &  \la_2  & &  \\ & \nearrow &  & \nwarrow & \\
\la_1\!\!\! &&&& \!\!\!\ka \\ & \nwarrow & & \nearrow& \\ & &
\ka_1 & &
\end{matrix}
\end{equation*}

\begin{Cla}\label{WWW}
$\big(\kova{\la_1}{\la_2}{\ka}\big)_{|\FF_{\ka}^{\bullet}} \in
\C\cdot\big(\au{\ka_1}{\la_1}\ad{\ka}{\ka_1}\big)_{|\FF_{\ka}^{\bullet}}
$.
\end{Cla}
\begin{proof}
We shall prove the equivalent statement:
\begin{equation*}
\au{\la_1}{\f}\kova{\la_1}{\la_2}{\ka} \in
\C\cdot\au{\la_1}{\f}\au{\ka_1}{\la_1}\ad{\ka}{\ka_1}
~\big(=\C\cdot\au{\ka_1}{\f}\ad{\ka}{\ka_1} \text{~by lemma
\ref{composition}}\big)
\end{equation*}
We begin by applying the r.h.s. to a cyclic element $\delta_{x_0}
\in \FF_{\ka}$:
\begin{equation*}\begin{split}
[\au{\ka_1}{\f}\ad{\ka}{\ka_1}\delta_{x_0}](y_0)&=\sum_{\substack{z
\subset y_0 \\
\ta(z)=\ka_1}}[\ad{\ka}{\ka_1}\delta_{x_0}](z)=\sum_{\substack{z
\subset y_0 \\ \ta(z)=\ka_1}}\sum_{\substack{x \supset z
\\ \ta(x)=\f}}\delta_{x_0}(x)\\
&=\sum_{\substack{z \subset y_0 \\
\ta(z)=\ka_1}}{\mathbf{1}}_{\{z|z \subset x_0, \ta(z)=\ka_1\}}
=\big|\{z|z \subset y_0 \cap x_0, \ta(z)=\ka_1\}\big|\\
&=\left\{
\begin{array}{ll}
    \big({\ka \atop \ka_1}\big) & \hbox{if $y_0 \supset x_0$} \\
    1 & \hbox{if $\ta(y_0 \cap x_0)=\ka_1$} \\
    0 & \hbox{otherwise} \\
\end{array}
\right.
\end{split}\end{equation*}
Let $w_0$, $x_0$ and $y_0$ be submodules of $F$ such that $w_0
\subset x_0 \cap y_0$ with types $\ta(w_0)=\ka_1$, $\ta(x_0)=\ka$
and $\ta(y_0)=\f$. Define:
\begin{equation*}
B=B_{x_0,y_0,w_0}=\{z | z \subset y_0, \ta(z)=\la_1,
\ta(z+x_0)=\la_2, z \cap x_0 = w_0 \} \qquad\qquad \qquad
\end{equation*}
By the second part of claim \ref{equivalence}:
\begin{equation*}
B=\{z |w_0 \subset z \subset y_0, \ta(z)=\la_1, \ta(\m z+ \m
x_0)=\m\la_2, \m z \cap \m x_0 = \m w_0\}
\end{equation*}
This description of $B$ together with the symmetricity of $y_0$
implies that $|B|$ depends only on the types of $\m x_0$ and
$w_0$, denote it by $b_{\m \ka,\ka_1}$. Applying the l.h.s. to
$\delta_{x_0}$ yields:
\begin{equation*}\begin{split}
[\au{\la_1}{f}\kova{\la_1}{\la_2}{\ka}\delta_{x_0}](y_0)&=\sum_{\substack{z
\subset y_0 \\ \ta(z)=\la_1}}
[\kova{\la_1}{\la_2}{\ka}\delta_{x_0}](z)=\sum_{\substack{z
\subset y_0 \\ \ta(z)=\la_1}}\sum_{\substack{x
\\ \ta(z+x)=\la_2}}\delta_{x_0}(x)\\&=\big|\{z| z \subset y_0, \ta(z)=\la_1, \ta(z+x_0)=\la_2
\}\big|\\&=\left\{
\begin{array}{ll}
    \Big({\ka \atop \ka_1}\Big)b_{\m \ka,\ka_1} & \hbox{if $y_0 \supset x_0$} \\
    b_{\m \ka,\ka_1} & \hbox{if $\ta(y_0 \cap x_0)=\ka_1$} \\
    0 & \hbox{otherwise} \\
\end{array}
\right.
\end{split}\end{equation*}

\noindent Hence, the two operators differ by a constant.

\end{proof}

\noindent \textbf{Step 3}

\noindent We should prove that the triple $\ka \le \e \le \f$ is
good when both $\e$ and $\f$ are symmetric types. This is
precisely the assertion of \cite[Theorem 6]{p1}. In particular
see relation (8) in the proof. The only delicate point which
deserves a remark, is the duality axiom which is used in the
proof and should be justified. Indeed, in a module of type $\ka$,
the number of submodules of type $\alpha$ equals to the number of
submodules of co-type $\alpha$ (cf. \cite[II.1]{macdonald}).

\subsection{Explicit calculation} Recall the notations of \S\ref{explicit.1}. For partitions $\la,\nu \in \cup\La_j$ let
\[ \langle\la,\nu\rangle=\sum{\la_i\nu_i}.\]

\begin{Cla}
\begin{align*}
(1) &\qquad \Bigl({\nu \atop \la}\Bigr)=
q^{\langle\nu'-\la',\la'\rangle}\prod_{i \ge 1}\Bigl[{\nu'_i-\la'_{i+1} \atop
\nu'_i-\la'_i}\Bigr].\\
(2) &\qquad
\cont{\la}{\nu}{\f}=q^{\langle\f'-\nu',\nu'-\la'\rangle}\Bigl[{\f_1'-\la_1'
\atop \f_1'-\nu_1'}\Bigr]\prod_{i \ge 1}\Bigl[{\nu'_i-\la'_{i+1}
\atop
\nu'_i-\nu'_{i+1}}\Bigr].\\
(3) &\qquad
\justavoid{\nu}{\la}_\Phi= q^{\langle\nu',{\la'_1}^\len\rangle}\Bigl({\len^{n-\la'_1} \atop
   \nu}\Bigr). \qquad \qquad \qquad \qquad \qquad
\qquad \qquad \\
(4) &\qquad
\avoid{\nu}{\f}{\la}_{\Phi}=q^{\langle\Phi'-\phi',\phi'-\nu'\rangle}\Bigl[{\Phi'_1-\nu'_1-\la'_1
\atop \phi'_1-\nu'_1}\Bigr]. \qquad \qquad \qquad \qquad \qquad
\qquad \qquad
\end{align*}

\end{Cla}

\begin{proof}
A basic quantity, which all other quantities are scaled to, is the
cardinality of $\Hom(\OO_{\la},\OO_\nu)$ which is
denoted and computed by
\[ \hom(\la,\nu)=q^{\langle\la',\nu'\rangle} \]
The subset of all injective morphisms will be denoted
$\Hom^{1-1}(\OO_{\la},\OO_{\nu})$ and we will use
\[ \hom^{1-1}(\la,\nu)=
   |\Hom^{1-1}(\OO_{\la},\OO_{\nu})|
\]
For a given type $\la$, $\hom^{1-1}(\la,\la)$ is computed in
\cite[II.1]{macdonald} where it is denoted $a_\la(q)$. A similar computation shows that
\[ \hom^{1-1}(\la,\nu) =\prod_{i \ge 1}
   \frac{[\nu'_i-\la'_{i+1}]!}{[\nu'_i-\la'_i]!}\hom(\la,\nu).
\]
Observe that the map
\[ \Hom^{1-1}(\OO_{\la},\OO_{\nu})
   \rightarrow \Gr(\la,\OO_{\nu}), \qquad
   \psi \mapsto \mathrm{Im}(\psi)
\]
is $\hom^{1-1}(\la,\la)$ to one, thus
\[ \Bigl({\nu \atop
   \la}\Bigr)=\frac{\hom^{1-1}(\la,\nu)}{\hom^{1-1}(\la,\la)},
\]
and (1) follows. Given an $\OO$-module $f$ of type $\f$, counting
in two ways the size of the set
\[ \{x,y<f~|~x<y,~\ta(x)=\la,~\ta(y)=\nu \} \]
gives
\[ \Bigl({\f \atop \la}\Bigr) \cont{\la}{\nu}{\f} =
   \Bigl({\f \atop \nu}\Bigr)  \Bigl({\nu \atop \la}\Bigr),
\]
which combined with (1) proves (2).

In order to prove (3) we need a little preparation. Let $E$ be a
module of type $\Phi$. Let $z<E$ be a fixed module of type $\la$.
Denote by $E_1$ a module of type $\len^{\la'_1}$ containing $z$.
Observe that $E_1$ is a direct summand of $E$. Fix a complimentary
direct summand to $E_2$ such that $E= E_1 \oplus E_2$. Denote the corresponding projections by $p_1$ and $p_2$. Assume
that a type $\nu$ is given. Let
\[
 X=\{x<E~|~\ta(x)=\nu,~x\cap z=0 \},
\]
in particular $\justavoid{\nu}{\la}=|X|$. Observe that the map
\[ \Hom(\OO_{\nu},E_1)\oplus
   \Hom^{1-1}(\OO_{\nu},E_2) \rightarrow X,
\]
\[   (\psi_1,\psi_2) \mapsto \mathrm{Im}(\psi_1+\psi_2)
\]
is $\hom^{1-1}(\nu,\nu)$ to one. This map is indeed into $X$ as
\[ \mathrm{Im}(\psi_1\oplus\psi_2) \cap z \subset
   \mathrm{Im}(\psi_1\oplus\psi_2) \cap E_1 \simeq \Ker (\psi_2) = (0).
\]
It is onto $X$ as for a given
$\psi\in\Hom^{1-1}(\OO_{\nu},E)$, with
$\mathrm{Im}(\psi)\in X$ we have $\Ker(p_2\circ \psi)\simeq
\mathrm{Im}(\psi)\cap E_1$ which is $(0)$, as the $\m$-torsion of
$E_1$ is equal to the $\m$-torsion of $z$. We get that
\[
\justavoid{\nu}{\la}= \Bigl({\len^{n-\la'_1} \atop
   \nu}\Bigr)\hom(\nu,\len^{\la'_1}),
\]
which proves (3). Counting in two ways the size of the set
\[ \{x,y<E~|~x<y,~y\cap z=0,~\ta(x)=\nu,~\ta(y)=\f \} \]
gives
\[ \justavoid{\f}{\la}_\Phi\Bigl({\f \atop \nu}\Bigr) =
   \justavoid{\nu}{\la}_\Phi\avoid{\nu}{\f}{\la}_\Phi
\]
which combined with (3) proves (4).

\end{proof}

\begin{Cor}
\[
A_{\la\ka}=q^{\langle\f'-\la',\Phi'-\phi'+\la'-\ka'\rangle}\Bigl[{\f_1'-\ka_1'
\atop \f_1'-\la_1'}\Bigr]\Bigl[{\Phi'_1-\la'_1-\ka'_1
\atop \phi'_1-\la'_1}\Bigr]\prod_{i \ge 1}\Bigl[{\la'_i-\ka'_{i+1}
\atop
\la'_i-\la'_{i+1}}\Bigr].
\]
\end{Cor}



\section{From the finite Grassmannians to $\Scal(\Gr(m,n,\Fb))$}\label{section:infinite}

Recall that $I_\len$ denotes the
kernel of the reduction of $\PGL(\OO)$ modulo $\m^\len$. The group $\PGL(\OO)$,
being the inverse limit of the finite groups $\gln$, enjoys the
property that each of its continuous irreducible complex
representations has a {\em level}. That is, the first nonnegative
integer $\len$ such that $I_{\len+1}$ acts trivially. It follows
that there exist a natural filtration
\[
(0) \subset \Scal(\Gr(m,n,\Fb))^{I_1} \subset \Scal(\Gr(m,n,\Fb))^{I_2}
\subset \cdots \subset \Scal(\Gr(m,n,\Fb)).
\]
The $\len$-th term in this filtration consists of all the
irreducible components of the representation which have level at
most $\len-1$, and thus can be regarded as a representation of
$\gln$. Denote this representation by $\bar{\rho}_{\len^m}$. Since
each irreducible constituent is contained in some
$\bar{\rho}_{\len^m}$, we get
\begin{equation}\label{limit.1}
\Scal(\Gr(m,n,\Fb)) \simeq \varinjlim \bar{\rho}_{\len^m}
\end{equation}
as $\GL_n(\OO)$ representations.

\begin{Lem}\label{equiv} $(\bar{\rho}_{\len^m},\Scal(\Gr(m,n,\Fb))^{I_\len}) \simeq (\rho_{\len^m},\FF_{\len^m})$.
\end{Lem}

\begin{proof} Follows from the natural identification
\[
I_\len \backslash \Gr(m,n,F) \simeq I_\len \backslash
\GL_n(F)/P_m(F) \simeq \GL_n(\OO_\len)/P_m(\OO_\len) \simeq
X_{\len^m},
\]
where $P_m$ is the appropriate parabolic group, and the
isomorphism
\[
\Scal(\Gr(m,n,\Fb))^{I_\len} \simeq \C(I_\len \backslash
\Gr(m,n,\Fb)).
\]
\end{proof}
Consequently, combining \eqref{limit.1} with Lemma \ref{equiv} proves Claim \ref{isom.grass.and.limit}. Note that the isomorphism in Lemma \ref{equiv} is algebraic, and
that there are two different inner products on $\rho_{\len^m}$
and on $\bar{\rho}_{\len^m}$, arising from the counting measure on
$X_{\len^m}$ or the projection of the (probability) Haar measure from $\GL_n(\OO)$ to $I_\len
\backslash \Gr(m,n,F)$, respectively. Throughout the finite
analysis we kept the former inner product, while at this stage of transferring the results to the infinite Grassmann representation, we
should keep track of the appropriate normalization.

The whole study included here hinges on the pro-finite nature of
the ring of integers $\OO$, and hence of all groups, spaces and
algebras defined over it, summarized by
\begin{align*}
&\mathrm{Groups}& &\PGL(\OO) \simeq \varprojlim \PGL(\Ok) \\
&\mathrm{Spaces}& &\Gr(m,n,\Fb) \simeq \varprojlim \Gr(m,n,\Ok)\\
&\mathrm{Representations}& &\Scal(\Gr(m,n,\Fb)) \simeq \varinjlim \FF_{\len^m}\\
&\mathrm{Algebras}& &\HH_m \simeq \varinjlim
\End_{\PGL(\OO)}(\FF_{\len^m})
\end{align*}
and explained in detail below.

\subsection{Lifting the finite spaces, algebras and functions}\label{gen}

Let $\pi_\len:X_{\len^m}\rightarrow
X_{(\len-1)^m}$ be the natural quotient maps. As $\GL_n(\OO)$-spaces
we have
\[
\Gr(m,n,\Fb) \simeq \varprojlim X_{\len^m}.
\]
Using the identification in Claim
\ref{ident.cosets.types}, the maps $\pi_\len$ descent to maps
\[
\pi_\len : \Lambda^\len_m \rightarrow \Lambda^{\len-1}_m ,\qquad
\tau(\OO_\la)\mapsto
   \tau(\OO_\la/\m^{\len-1}\OO_\la)
\]
which are easily described in transposed coordinates by
\[ (\lambda'_1,\ldots,\lambda'_{\len-1},\lambda'_\len) \mapsto
   (\lambda'_1,\ldots,\lambda'_{\len-1}). \]
Consider the set $\coprod_{\len\geq 0} \Lambda^\len_m$, and endow it
with a graph structure by connecting each $\la \in \Lambda^\len_m$
with its image $\pi_\len(\la) \in \Lambda^{\len-1}_m$. This graph is a
rooted tree, the root being the empty partition in $\Lambda^0_m$.
The inverse limit $\varprojlim\Lambda^\len_m$ can be identified with the space of ends
of this tree. The obvious sections $\Lambda^{\len-1}_m\rightarrow
\Lambda^\len_m$ given by $(\lambda'_1,\ldots,\lambda'_{\len-1}) \mapsto
   (\lambda'_1,\ldots,\lambda'_{\len-1},0)$ give at the limit an imbedding of
$\Lambda_m$ in $\varprojlim\Lambda^\len_m$. Thus, $\Lambda_m$ can be
identified with an open and dense subset of
$\varprojlim\Lambda^\len_m$. Let $\overline{\N}$ stand for the one
point compactification of $\N$. $\Lambda_m$ is naturally imbedded
in $\N^m$. We denote by $\overline{\Lambda}_m$ its closure in
$\overline{\N}^m$. It is easily seen that $\overline{\Lambda}_m$
can be identified with $\varprojlim \Lambda^\len_m$. We summarize
this discussion by

\begin{Pro}\label{spaces} $\Gr(m,n,\Fb)\times_{\GL_n(\OO)} \Gr(m,n,\Fb) \simeq \overline{\Lambda}_m \simeq \varprojlim \Lambda^\len_m$.
\end{Pro}
\begin{proof}
The only nontrivial issue left to address is the fact that the
first identification is also topological. The topology on the
l.h.s is the quotient topology. The quotient map from $\Gr(m,n,\Fb) \times
\Gr(m,n,\Fb)$ to $X_{\len^m} \times X_{\len^m}$ is continuous and
$\PGL(\OO)$-equivariant. The limit map becomes continuous and
well-defined on the quotient.
\end{proof}

\begin{Rem}\qquad
\begin{enumerate}
\item As topological spaces, $\overline{\Lambda}_m \setminus \Lambda_m
\simeq \overline{\Lambda}_{m-1}$, thus $\coprod_{i=0}^m \Lambda_i$
is a stratification of $\overline{\Lambda}_m$.

\item $\Gr(m,n,\Fb)$ carries a $\GL_n(\OO)$-invariant measure. Consequently, also does
$\overline{\Lambda}_m$, \\ and $\Lambda_m$ is of full measure
inside $\overline{\Lambda}_m$. This measure is computed in
\cite[\S2.2]{onn1}.
\end{enumerate}
\end{Rem}

\medskip

The maps $\pi_\len$ give rise to inclusions of the (finite
dimensional) spaces $i_\len: L^2(X_{(\len-1)^m})\rightarrow
L^2(X_{\len^m})$, where the notation $\FF(X)$ is replaced by $L^2(X)$
to emphasize that the inner product structure is induced from the
Haar measure, rather than the counting measure. The adjoint
transformation, $i_\len^*$, is the orthogonal projection on the
$I_{\len-1}$ invariants. In the limit we get the vector space of
Bruhat-Schwartz (=locally constant) functions:
\[
\Scal(\Gr(m,n,\Fb)) \simeq \varinjlim L^2(\Gr(m,n,\Fb))^{I_\len} \simeq \varinjlim
L^2(X_{\len^m}).
\]
$L^2(\Gr(m,n,\Fb))$ is the completion $\Scal(\Gr(m,n,\Fb))$, or alternatively, the direct limit in the category of Hilbert spaces. The inclusions $i_\len$
also give embeddings of the Hecke algebras
$\HH_{\len^m}=\End_{\gln}\left(L^2(X_{\len^m})\right)$ given by
\[
\HH_{(\len-1)^m} \rightarrow \HH_{\len^m},  \qquad h\mapsto i_\len \circ h
\circ i_\len^*.
\]
Recall by \eqref{identification.hecke} that as vector spaces $\HH_{\len^m} \simeq \FF(\Lambda^\len_m)$, and that under this isomorphism the operator
$\gb_\la \in \HH_{\len^m}$ ($\la \in \La_m^\len$) corresponds to the delta function $\delta^\len_{\la} \in \FF(\Lambda^\len_m)$ supported on $\la$.
\begin{Cla}
If $\la'_{\len-1}=0$ then $i_\len \circ \delta^{\len-1}_{\la} \circ
i_\len^*=\delta^\len_{\la}$.
\end{Cla}
\begin{proof}
The condition $\la'_{\len-1}=0$ implies that $\pi_\len^{-1}(\la)$ is the
singleton $\{\la\}\subset \La_m$, and the claim follows.
\end{proof}

It follows that the image of $\delta_{\la}^\len$ inside $\varinjlim \HH_{\len^m}$
stabilizes for $\len$ large enough, thus determining an element
$\delta_{\la}=\varinjlim \delta^\len_{\la}$,
where $\delta_{\la}$ is the delta function supported at $\la$,
viewing $\la$ as an element of $\bar{\La}_m$ via $\La_m \subset
\bar{\La}_m$. The identifications of Proposition \ref{spaces} give three ways to look at
\[
\HH_m=\End_{\GL_n(\OO)}\left(\Scal(\Gr(,m,n,\Fb))\right),
\]
namely,
\[
\HH_m \simeq \Scal\left(\Gr(m,n,\Fb)\times_{\GL_n(\OO)} \Gr(m,n,\Fb)\right) \simeq \Scal(\overline{\Lambda}_m) \simeq \varinjlim \HH_{\len^m}.
\]
Here $\Scal(\overline{\Lambda}_m)$ is the space of locally
constant functions on $\overline{\Lambda}_m$. The limit algebra
structure obviously coincides with the operator algebra structure
of $\Scal(\overline{\Lambda}_m)$, arising when viewing its
elements as convolution operators on $\Scal(\Gr(m,n,\Fb))$. Denote by
$\FF(\Lambda_m)$ the space of finitely supported functions on
$\Lambda_m$. As $\Lambda_m$ is discrete in $\overline{\Lambda}_m$,
$\FF(\Lambda_m)$ is imbedded in $\Scal(\overline{\Lambda}_m)$. As
$\Lambda_m$ is dense in $\overline{\Lambda}_m$, $\FF(\Lambda_m)$
is dense in $\Scal(\overline{\Lambda}_m)$ too. Consequently, the
algebraic structure of $\Scal(\overline{\Lambda}_m)$ is determined
by $\FF(\Lambda_m)=\Span\{\delta_{\la}:\la \in \Lambda_m\}$.

\subsection{Transition matrices}\label{explicitformula}
We are finally in a position to collect the pieces, and write down
explicitly the transition matrix between the delta functions basis
of $\HH_m$ and the idempotents of $\HH_m$. In order to do that we introduce an intermediate
basis, which is the limit of the (normalized image in $\HH_m$) of
the finite levels cellular bases. For $\la \in \La_m$ define
\[
\bar{\cb}_{\la}=\varinjlim_\len \frac{\cb_{\la}}{\bigl({\Phi \atop \phi}\bigr)},
\]
$\bar{\gb}_\la$ the operator which corresponds to $\delta_\la$, and $\bar{\eb}_\la$ the image of $\eb_\la$ in $\HH_m$. Combining the results of \S\ref{sec:c-g}-\ref{sec:c-e} with the above discussion gives
\begin{The}\label{main.theorem}
\begin{align*}
\bar{\gb}_{\la}&=\sum_{\ka \ge \la \ge \m\ka} \hat{\mu}(\la,\ka)
\bar{\cb}_{\ka} \\
 &= \sum_{\ka \ge \la \ge \m\ka} (-1)^{|\nu|-|\la|}
q^{n(\nu)-n(\la)} \prod \Bigl[ {\nu'_i-\nu'_{i+1} \atop
\nu'_i-\la'_i} \Bigr] \bar{\cb}_{\ka}\\ \\
\bar{\cb}_{\ka}&=\sum_{\nu \le \ka}
\frac{\cont{\nu}{\ka}{\f}\avoid{\ka}{\f}{\nu}_{\Phi}}{\bigl({\Phi
\atop \phi}\bigr)}\bar{\eb}_{\nu}\\ &=\sum_{\nu \le \ka}
q^{-(n-2m)|\ka|-m|\nu|-<\ka',\ka'-\nu'>} \frac{\Bigl[{ m-\nu'_1
\atop m-\ka'_1} \Bigr]\Bigl[ {n-\ka'_1-\nu'_1 \atop m-\ka'_1}
\Bigr]}{\Bigl[ {n \atop m} \Bigr]} \prod \biggr[
{\ka'_i-\nu'_{i+1} \atop \ka'_i-\nu'_i} \biggl] \bar{\eb}_{\nu}
\end{align*}
Which together give the desired $(\bar{\gb}-\bar{\eb})$ transition
matrix.
\end{The}


\section{Grassmann representation over nonsymmetric modules and open problems}\label{problems}

\subsection{Complexity of the representations $\FF_\la$} The focus of this paper is on the representations $\FF_\la$ with $\la=\len^m$. They enjoy the property
that their decomposition into irreducible constituents is of combinatorial nature, in particular, independent of
the ring $\OO$. For arbitrary (non-rectangular) types there is a strong dependence on the ring, and a highly
non trivial problem is

\begin{Prob} Decompose $\FF_\la$ into irreducible constituents for any $\la \in \La_m$.
\end{Prob}

The nontriviality of the problem is demonstrated in the next proposition. Let $B_\jmath$ denote the subgroup of upper triangular
matrices in $\GL_n(\OO_\jmath)$, that is, the stabilizer of a full flag of $\OO_\jmath$-free submodules in $\OO_\jmath^n$.

\begin{Pro} If all parts of $\la \in \La_n$ are pairwise unequal with smallest part $\la_n=\jmath$ and largest part $\la_1=\len$, then the $\gln$-representation $\FF_\la$ contains $\text{Ind}_{B_\jmath}^{G_{\jmath^n}}(1)$, where the action of $\gln$ on the latter is via reduction modulo $\m^\jmath$.
\end{Pro}

To get an idea of the complexity of
$\text{Ind}_{B_\jmath}^{G_{\jmath^n}}(1)$, the reader is referred
to \cite{CN} in which the case $n=3$ is studied. Though we know
very little about these arbitrary Grassmannians, they can be used
to pin down the irreducible representations which occur in the
Grassmann representation studied in the current paper: the essence
of the labeling in \eqref{decomposition} of \S\ref{introduction}
comes from the following theorem.

\begin{The}\cite{p1} There exist a family $\big\{\UU_{\la}^\Fb~|~\la \in
\La_m \big\}$ of irreducible representations of $\PGL(\OO)$ such that
\begin{enumerate}
    \item $\Scal(\Gr(m,n,\Fb)) = \bigoplus_{\la \in \Lambda_m}
\UU_{\la}^\Fb$.
    \item $\langle \UU_{\la}^\Fb,\FF_{\mu} \rangle =|\{\la \hookrightarrow \mu
    \}|$. I.e., the multiplicity of $\UU_{\la}^\Fb$ in $\FF_{\mu}$
    is the number of nonequivalent embeddings of a module
    of type $\la$ in a module of type $\mu$.
\end{enumerate}
In particular $\UU_\la^\Fb$ occurs both in $\Scal(\Gr(m,n,\Fb))$ and in $\FF_{\la}$ with multiplicity one, and does not occur in
$\FF_\mu$ for $\la \nleq \mu$.

\end{The}

\subsection{Dimensions of $\UU_\la^\Fb$}

The dimensions of the representations $\UU^{\Fb}_\la$ were computed in \cite{OS}, using a sophisticated and heavy computational machinery
of degenerations of certain generalized quantum dimensions formulae to actual dimensions of the $\UU^{\Fb}_\la$'s. They are given by the
following formula \cite[\S4.3]{OS}.
\begin{equation}\label{padic-dimension}
\textup{Dim}_{\mathbb{C}}(\UU^{\Fb}_{\lambda})=t^{-(n-2m+1)|\lambda|-2(\rho,\lambda)}\left[
{m \atop \partial\lambda' }\right]_t
\frac{(t^{n-\lambda'_1-\lambda'_2+2};t)_{\lambda'_1+\lambda'_2}}{(t^{m-\lambda'_1+1};t)_{\lambda'_1}}
\frac{(1-t^{n-2\lambda'_1+1})}{(1-t^{n+1})}
\end{equation}
for $\lambda\in\Lambda_n$.

Here $\partial \la'=(\la'_j-\la'_{j+1})_{j\ge0}$, $\rho=(n-1,n-2,\ldots,0)$, $(a;t)_j=\prod_{i=0}^{j-1}(1-t^{i}a)$ and $t=|\OO/\m|^{-1}$. We used $t$ instead of $q$ to match the notation of \cite{OS}.

\begin{Prob} Compute the dimensions of $\UU^{\Fb}_\la$ ($\la \in \La_m$) directly.
\end{Prob}

\subsection{Heisenberg-like relations on the complete lattice of submodules}
Let $\LL$ denote the lattice of submodules in $\OO_\len^n$ and let $\La=\La_n^\len$. Let $\FF(\LL)=\oplus_{\la \le \len^n} \FF_\la$ stand for complex valued functions on $\LL$. For $x,y \in \LL$ we use the notation $y \gtrdot x$
whenever $y$ covers $x$ (i.e. $y/x$ is simple). Define the
following 'lowering' and 'raising' operators on $\FF(\LL)$:
\begin{align*}
  \Db^{\flat}:\FF(\LL) \rightarrow \FF(\LL)& &\Db^{\sharp}:\FF(\LL) \rightarrow \FF(\LL)\\
  \Db^{\flat}f(x)=\sum_{y \gtrdot x}f(y) &      & \Db^{\sharp}f(x)=\sum_{y \lessdot x}f(y).
\end{align*}
Observe that $\Db^{\flat}$ and $\Db^{\sharp}$ are adjoints.
Indeed, this follows from $\Db^{\flat}=\sum_{\mu \lessdot \la}
\ad{\la}{\mu}$ and $\Db^{\sharp}=\sum_{\mu \gtrdot \la}
\au{\la}{\mu}$.

\begin{Pro}\label{Heisenberg}
 $(\Db^{\flat}\Db^{\sharp}-\Db^{\sharp}\Db^{\flat})_{|\FF_{\la}}=b_{\la} \cdot \mathbf{Id}_{\FF_{\la}}$
\end{Pro}

\begin{proof} Using the definition of $\Db^{\flat}$ and $\Db^{\sharp}$, we need to show that
\begin{equation*}
 \sum_{y \gtrdot x} \sum_{z \lessdot y}{h(z)}-\sum_{y \lessdot
 x}\sum_{z \gtrdot y}{h(z)} = b_{\la} \cdot h(x) \qquad \forall
 x \in X_{\la}
\end{equation*}

For any subset $\Sigma \subset \La$ let $\sharp(\Sigma)$ be the
set of types which covers types from $\Sigma$ and $\flat(\Sigma)$
the set of types which are covered by types from $\Sigma$. For
$\la \in \La$ let
$\natural(\la)=\flat\sharp(\la)=\sharp\flat(\la)$.

First, we note that any $y \neq x$ from $\ta^{-1}(\natural\la)$,
appears exactly once in each of the summands on the left hand
side. Indeed, there is exactly one submodule (their join) which
covers both of them, and exactly one submodule (their meet) which
is covered by both of them. Hence such pairs do not contribute to
the left hand side.

Second, $y=x \in \ta^{-1}(\natural\la)$ appears $u_{\la}$ times in
the first summand, and $l_{\la}$ times in the second summand,
where for a fixed $z_0 \in \ta^{-1}(\la)$:
\begin{align*}
u_{\la}= |\{y| y \gtrdot z_0 \}| &\qquad \qquad \qquad  l_{\la}=
|\{y| y \lessdot z_0\}|
\end{align*}

\noindent It follows that $b_{\la}=u_{\la}-l_{\la}$.
\end{proof}

The scalars $b_\la$ can be easily computed. They are given by
\begin{equation*}
b_{\la}=\frac{q^{n-\rk\la}-q^{\rk\la}}{q-1}
\end{equation*}
which follows from
\begin{align*}
l_{\la}&=|\{y| y \lessdot
z_0\}|=\big|\mathbb{P}_{\kk}^{\rk(\la)-1}\big|=\frac{q^{{\rk(\la)}}-1}{q-1}\\
u_{\la}&=|\{y| y \gtrdot z_0
\}|=\big|\mathbb{P}_{\kk}^{n-\rk(\la)-1}\big|=\frac{q^{n-\rk(\la)}-1}{q-1}.
\end{align*}

\subsection{More questions regarding the Hecke algebra $\HH_{\len^m}$}

The transition matricx (\textbf{c-e}) is given explicitly by a combinatorial data (Theorem \ref{cell-idem}).
Examples imply that this should also be the case for the
transition matrix (\textbf{e-c}), and it would be
interesting to find such interpretation.

\begin{Prob} Invert the relation (\textbf{c-e}).
\end{Prob}

As mentioned in the introduction the case $\len=1$ is well studied \cite{dunkl,delsarte-H}. In {\em loco citato} the set $X_{1^m}$ is studied as an association scheme
(the $q$-Johnson scheme). The graph structure is defined by: two points $x,y \in X_{1^m}$ are connected with an edge if $x \cap y$ is of codimension one in each of them. The Laplacian $\Delta$ on this graph, defined by $\Delta h(x)=\sum_{x \sim y}h(y)$ (for $h \in \FF_{1^m}$), is nothing but the operator $\gb_{1^{m-1}}$. It turns out that $\Delta$ generates the Hecke algebra, i.e. $\HH_{1^m}=\C[\Delta]$. The following is a conjectural generalization of this fact.

\begin{Prob} Prove that $\{\gb_{(\len^{m-1},j^1)}~|~ j=0,\ldots,\len-1\}$ generate $\HH_{\len^m}$ as an algebra.
\end{Prob}

\subsection{Other algebraic groups}
The case $\len=1$ admits generalizations to other algebraic
groups, see e.g. \cite{SD1}. It would be interesting to generalize
these further to other classical groups, for example

\begin{Prob} Study the natural Grassmann representation of the symplectic group arising from its action on Lagrangian subspaces of $\len^{2n}$.
\end{Prob}


\appendix
\section{Local rings and discrete valuation rings}\label{dvr}
In this appendix we prove some claims regarding local rings. All
modules under consideration are assumed to be of finite rank. Let $R$ be a local ring with a maximal ideal $\m$
and let $x, y$ be $R$-modules. Let $\bar{x}=x/\m x$. By Nakayama's
lemma $x \rightarrow y$ is onto if and only if the induced map
$\bar{x} \rightarrow \bar{y}$ is onto. Equivalently $\rk(x)=\dim
x/\m x$.

\begin{Cla}
Let $R$ be a local ring with maximal ideal $\m$. Let $z$ be an
$R$-module and $x$, $y$ two submodules of $z$. Then:
\begin{equation*}
\rk(x+y)=\rk(x)+\rk(y)-\rk(x \cap y)+\dim\Big(\frac{\m x \cap \m
y}{\m(x \cap y)}\Big)
\end{equation*}
\end{Cla}

\begin{proof} Using the equalities $\m(x+y)=\m x + \m y$ and $\m(x \oplus y)=\m x \oplus \m
y$, we obtain a commutative diagram with exact columns and rows
\eqref{homolog} and exact sequence \eqref{exact}:
\begin{equation}\label{homolog}
\begin{matrix}
    &      &   0    &    & 0      &  &   0 & & \\
    &      &   \dar &    & \dar   &  & \dar & &\\
 0  & \rar &  \m x\cap\m y & \rar &  \m x\oplus \m y &\rar &\m x
 +\m y & \rar &0 \\
    &      &   \dar &    & \dar   &  & \dar & &\\
 0  & \rar &   x \cap y & \rar &  x \oplus y &\rar & x + y & \rar &0 \\
    &      &   \dar &    & \dar   &  & \dar & &\\
 0  & \rar &  \frac{x \cap y}{\m x\cap\m y} & \dashrightarrow &  \frac{x \oplus y}{\m x \oplus \m y}
 &\dashrightarrow &\frac{x+y}{\m
 (x + y)} & \rar &0 \\
    &      &   \dar &    & \dar   &  & \dar & &\\
    &      &   0    &    & 0      &  &   0 & & \\
\end{matrix}
\end{equation}
(the exactness of the dashed row follows from the obvious
exactness of the other rows and columns)

\begin{equation}\label{exact}
 0  \rar \frac{\m x \cap \m y}{\m (x \cap y)} \rar \frac{x \cap y}{\m (x \cap y)} \rar \frac{x \cap y}{\m
 x \cap \m y}  \rar  0
\end{equation}
and obtain:
\begin{equation*}\begin{split}
\rk(x)+\rk(y)=\rk(x \oplus y)
&\overset{\eqref{homolog}}{=}\rk(x+y)+\dim\frac{x \cap y}{\m x \cap
\m y} \\ &\overset{\eqref{exact}}{=}\rk(x+y)+\rk(x \cap
y)-\dim\Big(\frac{\m x \cap \m y}{\m(x \cap y)}\Big)
\end{split}\end{equation*}
\end{proof}

\begin{Cla}\label{exact1}
Let $R$ be a local ring with maximal ideal $\m$. Let $z$ be an
$R$-module and $x$, $y$ two submodules of $z$. Then:
\begin{equation*}
\rk(x)=\rk(x+y) \Longleftrightarrow \rk(x \cap y)=\rk(y)
\mathrm{~and~} \m x \cap \m y = \m (x \cap y)
\end{equation*}
\end{Cla}
\begin{proof} ($\Leftarrow$) clear. ($\Rightarrow$) if
$\rk(x)=\rk(x+y)$ we have that the two non-negative terms
$\rk(y)-\rk(x \cap y)$ and $\dim\big(\frac{\m x \cap \m y}{\m(x \cap
y)}\big)$ must sum up to zero.
\end{proof}

We now specialize to the situation to which we apply these
assertions.
\begin{Cla}\label{equivalence}
Let $\OO$ be a discrete valuation ring with maximal ideal $\m$. Let
$z$ be a finite $\OO$-module and $x$, $y$ two submodules of $z$.
Assume that $\rk(x+y)=\rk(x)$ and $x+y$ covers $x$.
\begin{enumerate}
\item  $\ta(x \cap y)$ depends only on $\ta(x)$, $\ta(y)$
and $\ta(x+y)$.

\item  $\ta(x+y)=\la \Longleftrightarrow \ta(\m x +\m y)=
\m \la$, $\m x \cap \m y = \m (x \cap y)$ and $\m(x \cap y)
\subsetneqq \m y$.

\end{enumerate}
\end{Cla}

\begin{proof}\quad
One immediately verifies that for a module $w$, $\ta(w) =\la
\Leftrightarrow \dim \bar{w}=\la_1=\rk(\la)$ and $\ta(\m w)=\m \la$.
\begin{enumerate}
\item There exist a unique $i$ such that
$\m^{i}(x+y)/\m^{i}x \ne (0)$. By the isomorphism
$\m^{i}(x+y)/\m^{i}x \simeq \m^{i}y/\m^{i}(x \cap y)$ it is also the
unique $i$ for which $\m^{i}y/\m^{i}(x~\!\cap~\!y) \ne (0)$.

\item ($\Rightarrow$) Assume $\ta(x+y)=\la$. Clearly
$\ta(\m x +\m y)= \m \la$, and by Claim~\ref{exact1} also $\m x \cap
\m y = \m (x \cap y)$. Since $x \cap y \subsetneqq y$ and have the
same rank, the dimensions of $\overline{x \cap y}$ and $\bar{y}$
must be the same and $\m(x \cap y) \subsetneqq \m y$.

\noindent ($\Leftarrow$) The data on the right together with claim
\ref{exact1} implies that $\ta(\m(x+y))=\m\la$ and
$\rk(x+y)=\rk(x)$. This implies that $\ta(x+y)=\la$.

\end{enumerate}
\end{proof}

\vspace{\bigskipamount}

\begin{footnotesize}
\begin{quote}

Uri Bader\\
Department of Mathematics, Technion, Haifa 32000, Israel\\
{\tt uri.bader@gmail.com} \\

Uri Onn\\
Department of Mathematics, Ben-Gurion University of the Negev,\\
Beer-Sheva 84105, Israel\\
{\tt urionn@math.bgu.ac.il}

\end{quote}
\end{footnotesize}

\end{document}